%% file: arXiv_main.tex
\documentclass{revtex4-1}
\pdfoutput=1

\input{preamble}

\newcommand\bt{}
\newcommand\et{}
\newcommand\bp{}
\newcommand\ep{ Q.E.D. }

\begin{document}

\title{Ratio of Symmetries Between any two $n$-Node Graphs}


\author{Justus Isaiah Hibshman}
\email{jhibshma@nd.edu}
\affiliation{University of Notre Dame}

\begin{abstract}
Given any two graphs on the same vertex set, $G_1 = (V, E_1)$ and $G_2 = (V, E_2)$, along with the difference between the two graphs $\Delta = (E_1 \setminus E_2) \cup (E_2 \setminus E_1)$, we prove that the ratio of the sizes of the two graphs' automorphism groups is equivalent to the ratio of the sizes of $\Delta$'s automorphism orbits in $G_1$ and $G_2$ respectively. This result provides a link between graphs' symmetries that might otherwise seem to be unrelated.
\end{abstract}

\maketitle

\input{body}

\bibliographystyle{plain}  
\bibliography{references}

\newpage
\input{appendix}

\end{document}

%% file: preamble.tex
\usepackage{amsmath}
\usepackage{amssymb}

\input{basic_tikz_figures}
\usepackage{bbm} 

\newcommand*{\ie}{\textit{i.e.} }



%% file: basic_tikz_figures.tex
\usepackage{tikz,amsmath, amssymb,bm,color}
\usetikzlibrary{shapes,arrows}
\usetikzlibrary{calc}

\usetikzlibrary{matrix}
\usetikzlibrary{shapes.multipart}

\tikzstyle{basicnode} = [circle, thick, fill=lightgray!30, draw=black, text=black]
\tikzstyle{basicedge} = [draw=black, thin, text=black]
\tikzstyle{textnode} = [font=\small, text=black, draw=none]

\tikzstyle{highlightededge} = [draw=black, dashed, thick]

\tikzstyle{bluenode} = [circle, thick, fill=blue!30, draw=blue, text=black]

\tikzstyle{blackedge} = [draw=black, thick, text=black]
\tikzstyle{blueedge} = [draw=blue, thick, dashed, text=black]
\tikzstyle{rededge} = [draw=red, thick, dashed, text=black]

\tikzstyle{blackdiredge} = [->, draw=black, thick, text=black]
\tikzstyle{bluediredge} = [->, draw=blue, thick, dashed, text=black]
\tikzstyle{reddiredge} = [->, draw=red, thick, dashed, text=black]

\tikzstyle{blanknode} = [draw=none]

\usetikzlibrary {quotes}

%% file: body.tex
\section{Introduction}

Slight differences between two graphs' edge sets can entail slight or immense differences between the two graphs' automorphism groups.

For instance, the automorphism groups of a graph $G$ and of the same graph with one of its edges $e$ removed can be arbitrarily different. Hartke, Kolb, Nishikawa, and Stolee show this when they proved that given any two finite groups $A$ and $B$, which can be arbitrarily different, they can construct a graph $G$ which has some subgraph $G-e$, whose automorphism groups are isomorphic to $A$ and $B$ respectively~\cite{hartke2009automorphism}. This  holds for \emph{all} pairs of groups. In short, this previous result means that merely knowing a graph's automorphism group (the ``$A$'') tells you nothing about the automorphism group if an edge (or set of edges) is deleted (the ``$B$'').

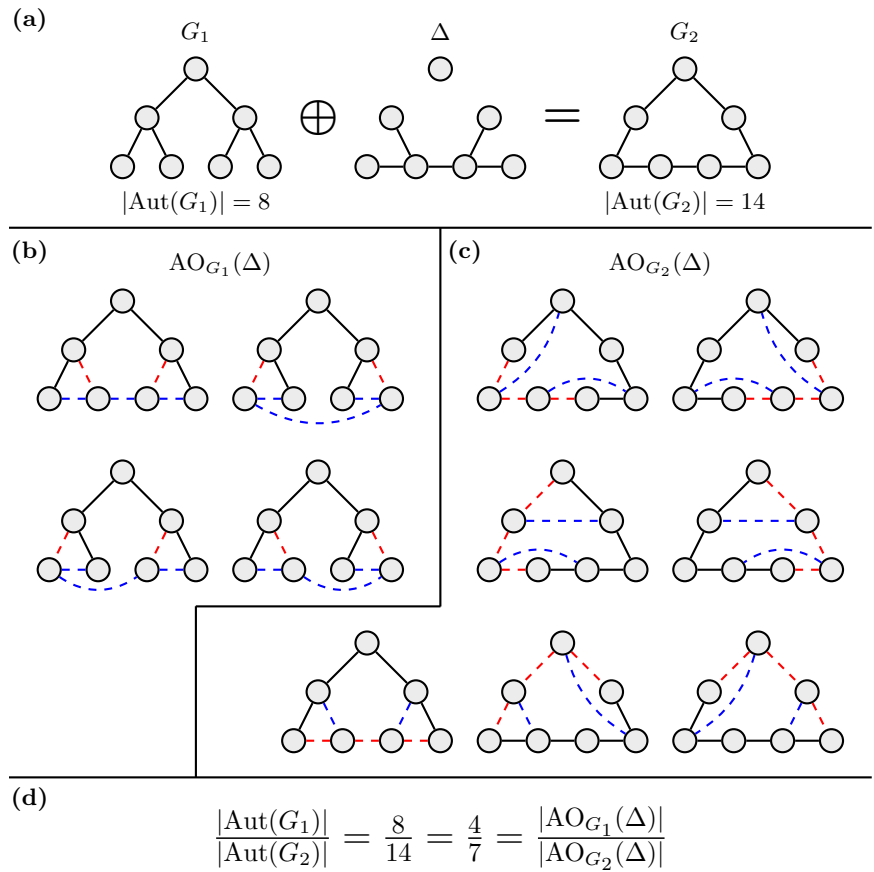
\begin{figure}[t]
\centering
\input{figures/main_result_01}
\caption{\textbf{Example -- Tree and Cycle} -- \textbf{(a)} shows two graphs $G_1$ and $G_2$, along with the difference between them $\Delta$. \textbf{(b) and (c)} show the automorphism orbits of $\Delta$ in $G_1$ and $G_2$ respectively. Red lines indicate edges (to be deleted). Blue lines indicate non-edges (to be added).
See Section~\ref{sec:prelim} and Figure~\ref{fig:example_edge_set_orbit} for more information on automorphism orbits of (non-)edge sets.
Finally, \textbf{(d)} shows our formula applied to this example.}\label{fig:main_example_01}
\end{figure}

In this paper, we look at any two graphs on the same node set and provide a link between the symmetries of the two graphs. Consider two graphs $G_1 = (V, E_1)$ and $G_2 = (V, E_2)$. We can note the difference between the graphs $\Delta = (E_1 \setminus E_2) \cup (E_2 \setminus E_1)$ and think of $\Delta$ as a way to transform $G_1$ into $G_2$ and vice versa by flipping every (non$-$)edge in $\Delta$. Now consider all the ways to transform $G_1$ into something isomorphic to $G_2$ where the transformation is structurally (\ie automorphically) equivalent to $\Delta$ in $G_1$; that set of transformations is the automorphism orbit of $\Delta$ in $G_1$. Likewise, consider the corresponding set for transforming $G_2$ into something isomorphic to $G_1$. In this paper we prove that the ratio of the sizes of the automorphism groups of $G_1$ and $G_2$ is equivalent to the ratio of the sizes of the automorphism orbits of $\Delta$ in $G_1$ and $G_2$ (i.e. the ratio of the sizes of the transformation sets). Our result is most easily understood and appreciated via a visual example, so we recommend beginning with Fig.~\ref{fig:main_example_01}.

The example illustrated in Fig.~\ref{fig:main_example_01}(\textbf{a}) shows a binary tree and its transformation into into a cycle via $\Delta$. \ref{fig:main_example_01}(\textbf{b}) illustrates all of the transformations from a tree to a cycle that are equivalent to $\Delta$. Likewise, \ref{fig:main_example_01}(\textbf{c}) shows all the equivalent transformations from the cycle to a tree. In this particular example, \ref{fig:main_example_01}(\textbf{d}) shows that the ratios of automorphism-group sizes is equal to the ratio of the sizes of the transformation sets. Indeed, in the present work, we show that
$\frac{|\text{Aut}(G_1)|}{|\text{Aut}(G_2)|} = \frac{|\text{AO}_{G_1}(\Delta)|}{|\text{AO}_{G_2}(\Delta)|}$
for any two graphs of the same number of vertices.

Note that the difference between the two graphs plays a special role in relating the graphs' symmetries. If we replaced $\Delta$ with some other edge set $S$, we would not necessarily get a link between the graphs' symmetries. For example, if $S$ consisted of the top two edges of the binary tree in Fig.~\ref{fig:main_example_01}, we would get that $\frac{|\text{AO}_{G_1}(S)|}{|\text{AO}_{G_2}(S)|} = \frac{1}{7}$, which does not equal the special $\frac{4}{7}$ symmetry ratio depicted in Fig.~\ref{fig:main_example_01} that we get from considering $\Delta$.

Our result applies to both directed and undirected graphs. We include an example figure of directed graphs in the appendix.

We suspect that this result may be relevant to two longstanding problems in graph theory: Graphical Enumeration~\cite{harary2014graphical} and the Graph Reconstruction Conjecture~\cite{bondy1977graph}.

Graphical enumeration is the task of calculating how many isomorphically distinct graphs exist on $n$ nodes. The calculation is equivalent to calculating how much automorphic symmetry the graphs on $n$ nodes have~\cite{harary1967number, harary2014graphical}, and our result provides insight into to how symmetry changes as a graph changes (\ie as edges are added). At present the problem is computationally intractable for even small $n$.

The graph reconstruction conjecture (and the related edge reconstruction conjecture) asks whether one can re-construct an $n$-node, $m$-edge graph after being shown the ``deck'' of $n$ (or $m$) graphs that arise from removing the graph's distinct vertices (or edges)~\cite{bondy1977graph, greenwell1969reconstructing, bollobas1990almost, kocay1987family, mcmullen2006graph, andersen1996set}. Removing a vertex from a graph involves removing the set of edges incident on that vertex. Since our result connects the symmetries of graphs before and after edges are removed, it may be useful in determining whether the symmetries of the deck can entail that a graph is reconstructible.

This paper next provides a preliminary introduction to the key concepts in Sec.~\ref{sec:prelim}. 
The main result is given in Sec.~\ref{sec:main_result} and the proof is offered in Sec.~\ref{sec:proof}.

\section{Preliminaries}\label{sec:prelim}

\subsection{Notation}\label{sec:shorthand_notation}

Let $S$ be a set and let $f : S \rightarrow S$ be a bijection. Given two elements $a, b \in S$, we use $f((a, b))$ to denote $(f(a), f(b))$.
Similarly, given a set of pairs $X \subseteq S \times S$ we use $f(X)$ to denote $\{f((c, d))\ |\ (c, d) \in X\}$.

Given two sets $A$ and $B$ we use $A \oplus B$ to denote $(A \setminus B) \cup (B \setminus A)$\footnote{Here $\oplus$ is analogous to the XOR function.}. Similarly, given a graph $G = (V, E)$ and set of node-pairs $P$, we use $G \oplus P$ to denote the graph $(V, E \oplus P)$.
\subsection{Iso- and Auto-morphisms}

Given a graph $G = (V, E)$ and a graph $G' = (V', E')$, an isomorphism between $G$ and $G'$ is a bijection $f : V \rightarrow V'$ such that $(a, b) \in E \leftrightarrow f((a, b)) \in E'$. Graphs $G$ and $G'$ are denoted to be isomorphic by writing $G \cong G'$.

An automorphism of $G$ is simply an isomorphism from $G$ to itself.

Let Aut($G$) denote the set of all automorphisms of $G$.

If a graph $G = (V, E)$ has colored edges, then we require that an automorphism only map edges of the same color to each other. Formally, if we have a set of colors $\mathcal{C}$ and edge coloring function $c: E \rightarrow \mathcal{C}$, then an automorphism of $G$ is a bijection $f: V \rightarrow V$ such that: $e \in E \rightarrow \left( f(e) \in E \land c(e) = c(f(e)) \right)$

\subsection{Automorphism Orbits}

Let $x$ be a graph \textit{entity} where $x$ could be a node ($x \in V$), a node pair ($x \in V \times V$), or a set of node pairs ($x \subseteq V \times V$); $x$ does not need to be an induced subgraph. The automorphism orbit of $x$ is the set of all entities that play the same structural role in $G$ that $x$ does:

$$\text{AO}_G(x) = \{f(x)\ |\ f \in \text{Aut}(G)\}$$

\begin{figure}[t]
\centering
\input{figures/edge_set_orbit_example}
\caption{\textbf{Example of the Automorphism Orbit and Stabilizer of a Set of Edges:} On the right, each distinct set of edges in the automorphism orbit is shown in a unique color and labeled with a letter. The multi-edges on the right are shown solely to indicate that the single edges on the left participate in multiple edge sets in the orbit. Note that sets such as $\{(1, 5), (5, 6)\}$ are not part of the example orbit. At the bottom, we can see the two automorphisms in the stabilizer listed; note that the only difference between the two automorphisms is that they swap the positions of nodes $3$ and $4$; everything else is constrained by the stabilized edge set.}\label{fig:example_edge_set_orbit}
\end{figure}
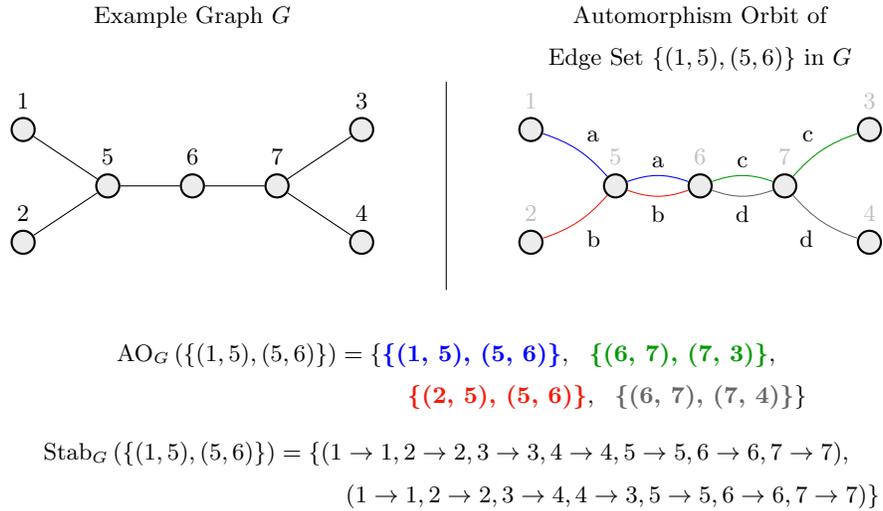

Note that if $x = (a, b)$, $x$ could be an edge in $E$ or $x$ could just as easily be a non-edge. Likewise, if $x$ is a set of node pairs, $x$ could be a set of edges, a set of non-edges, or a mixture. Lastly, note that the automorphism orbit of an entity $x$ contains the same kind of entities that $x$ is, per our notation in Section~\ref{sec:shorthand_notation}; the orbit of a node is a set of nodes, the orbit of an edge is a set of edges, etc. To see an example of the automorphism orbit of a set of edges, see Figure~\ref{fig:example_edge_set_orbit}.

\subsection{Stabilizers}

Graph automorphisms are realignments of a graph with itself; if for some graph entity $x$ in a graph $G$, an automorphism of $G$ realigns that $x$ with itself, then that automorphism is a stabilizer for $x$ in $G$~\cite{rose2009action}. Formally, the stabilizer set of an entity $x$ in a graph $G$ is the set:

$$\text{Stab}_G(x) = \{f\ |\ f \in \text{Aut}(G) \land f(x) = x\}$$

As with automorphism orbits, we can define stabilizers for nodes, edges, sets of edges, etc. To see an example of the stabilizer of a set of edges, see Figure~\ref{fig:example_edge_set_orbit}.

\section{Main Result}\label{sec:main_result}

\bt\label{thm:main_result}
For any two graphs $G_1 = (V, E_1)$, $G_2 = (V, E_2)$, and $\Delta=E_1 \oplus E_2$ it holds that:
\et
$$\frac{|\text{Aut}(G_1)|}{|\text{Aut}(G_2)|} = \frac{|\text{AO}_{G_1}(\Delta)|}{|\text{AO}_{G_2}(\Delta)|}$$

Remember that $\Delta$ is simply a set of node pairs and therefore every element in $\Delta$ can be thought of as an edge (or non-edge) in $G_1$ which is conversely a non-edge (or edge) in $G_2$. The above formula states that the ratio of automorphisms of $G_1$ to automorphisms of $G_2$ is equal to the ratio of the number of edge/non-edge sets that $\Delta$ is equivalent to in $G_1$ vs. in $G_2$.

Again, Figures~\ref{fig:main_example_01} and~\ref{fig:main_example_03} illustrate the result for undirected and directed graphs.

\section{Proof of Main Result}\label{sec:proof}

\bp
We begin our proof with a reminder concerning the orbit-stabilizer theorem~\cite{rose2009action}.

\ 

\textbf{Orbit-Stabilizer Theorem Applied to Graphs:} \emph{For any directed or undirected graph $G = (V, E)$, and for any node, node-pair, or set of node pairs $x$ in $G$:} $|\text{Aut}(G)| = |\text{AO}_G(x)| \cdot |\text{Stab}_G(x)|$




\ 

Let $G_1 = (V, E_1)$ and $G_2 = (V, E_2)$ be a pair of directed or undirected graphs. Let $\Delta = E_1 \oplus E_2$. Consider constructing a third graph $G_3 = (V, E_1 \cup E_2)$ wherein we color the edges as follows: $e \in E_1 \cap E_2$ means that $e$ is colored black; $e \in E_1 \setminus E_2$ means that $e$ is colored red; $e \in E_2 \setminus E_1$ means that $e$ is colored blue.

Note that $E_1 \setminus \Delta = E_2 \setminus \Delta = $ the set of black edges in $E_3$, $E_1 \cap \Delta = \Delta \setminus E_2 = $ the set of red edges in $E_3$, and $\Delta \setminus E_1 = E_2 \cap \Delta = $ the set of blue edges in $E_3$.

Now, we claim that $\text{Stab}_{G_1}(\Delta) = \text{Aut}(G_3) = \text{Stab}_{G_2}(\Delta)$. This can be seen if we expand the first two terms into their definitions:

\begin{equation}
\begin{aligned}
\text{Stab}_{G_1}(\Delta) = \{f\ |\ f \in & \text{Aut}(G_1) \land f(\Delta) = \Delta\} \\
= \{f\ |\ f:\ & V \rightarrow V \text{ is a bijection s.t.} \\ 
		& \left((a, b) \in E_1 \rightarrow (f(a), f(b)) \in E_1 \right)\\
		\land & \left((a, b) \in \Delta \rightarrow (f(a), f(b)) \in \Delta \right) \} \\
= \{f\ |\ f:\ & V \rightarrow V \text{ is a bijection s.t.} \\ 
		& \left((a, b) \in E_1 \cap \Delta \rightarrow (f(a), f(b)) \in E_1 \cap \Delta \right)\\
		\land & \left((a, b) \in E_1 \setminus \Delta \rightarrow (f(a), f(b)) \in E_1 \setminus \Delta \right) \\
		\land & \left((a, b) \in \Delta \setminus E_1 \rightarrow (f(a), f(b)) \in \Delta \setminus E_1 \right) \} \\
\ \ \ \ \ \ \ \ \ \ = \{f\ |\ f: & V \rightarrow V \text{ is a bijection s.t.} \\ 
		& \left((a, b) \in E_3 \text{ and red } \rightarrow (f(a), f(b)) \in E_3 \text{ and red} \right)\\
		\land & \left((a, b) \in E_3 \text{ and black } \rightarrow (f(a), f(b)) \in E_3 \text{ and black} \right) \\
		\land & \left((a, b) \in E_3 \text{ and blue } \rightarrow (f(a), f(b)) \in E_3 \text{ and blue} \right) \} \\
		= \text{Aut}(G_3&)
\end{aligned}
\end{equation}

The fact that $\text{Aut}(G_3) = \text{Stab}_{G_2}(\Delta)$ follows automatically from the symmetry of $G_1$ and $G_2$'s definitions relative to $G_3$. For visual intuition about this equality, observe our example figures~\ref{fig:main_example_01} and~\ref{fig:main_example_03} wherein $G_3$ can be seen drawn as the top-left element of $\text{AO}_{G_1}(\Delta)$.

Given that $\text{Stab}_{G_1}(\Delta) = \text{Stab}_{G_2}(\Delta)$, it follows from the orbit stabilizer theorem that:

\begin{equation}
\frac{|\text{Aut}(G_1)|}{|\text{AO}_{G_1}(\Delta)|} = |\text{Stab}_{G_1}(\Delta)| = |\text{Stab}_{G_2}(\Delta)| = \frac{|\text{Aut}(G_2)|}{|\text{AO}_{G_2}(\Delta)|}
\end{equation}

Re-arranging gives us our result:

\begin{equation}
\frac{|\text{Aut}(G_1)|}{|\text{Aut}(G_2)|} = \frac{|\text{AO}_{G_1}(\Delta)|}{|\text{AO}_{G_2}(\Delta)|}
\end{equation}

\ep


\section{Conclusion}\label{sec:conclusion}

We have provided a simple and concise proof of the fact that any two graphs' symmetries are related by the number of equivalent ways to transform the one graph into an isomorphic copy of the other (and vice versa).

In addition to the result's inherent beauty, we speculate that the result might prove useful for the surprisingly difficult task of calculating the total number of isomorphically distinct graphs on $n$ vertices, and it may also yield insights into the graph reconstruction conjecture, because in both problems, the relative symmetries of slightly modified graphs are deeply relevant.

\section*{Acknowledgements}
This research was supported by a grant from the US National Science Foundation (\#1652492).

Thanks to Tim Weninger and Daniel Gonzalez Cedre for their useful feedback in improving the presentation of this document. 


%% file: figures/main_result_01.tex
\begin{tikzpicture}[scale=0.65]


\node [textnode] at (0.6,0.25) {\textbf{(a)}};
\node [textnode] at (0.6,-4.5) {\textbf{(b)}};
\node [textnode] at (9.5,-4.5) {\textbf{(c)}};
\node [textnode] at (0.6,-15.7) {\textbf{(d)}};

\node [textnode] at (4,0) {$G_1$};
\node [textnode] at (6.5,-1.75) {\LARGE $\bm{\oplus}$};
\node [textnode] at (9,0) {$\Delta$};
\node [textnode] at (11.5,-1.75) {\LARGE $\bm{=}$};
\node [textnode] at (14,0) {$G_2$};

\begin{scope}[shift={(0,-3.5)}]
	\node [textnode] at (4, 0) {$|\text{Aut}(G_1)| = 8$};
	\node [textnode] at (14, 0) {$|\text{Aut}(G_2)| = 14$};
\end{scope}

\begin{scope}[shift={(0,-4.75)}]
	\node [textnode] at (4.5, 0) {$\text{AO}_{G_1}(\Delta)$};
	\node [textnode] at (13.5, 0) {$\text{AO}_{G_2}(\Delta)$};
\end{scope}

\node [textnode] at (9,-16.5) {\Large $\frac{|\text{Aut}(G_1)|}{|\text{Aut}(G_2)|} = \frac{8}{14} = \frac{4}{7} = \frac{|\text{AO}_{G_1}(\Delta)|}{|\text{AO}_{G_2}(\Delta)|}$};

\begin{scope}[shift={(4,-0.75)}]
 	\node [basicnode] (v1) at (0,0) {};
	\node [basicnode] (v2) at (-1,-1) {};
	\node [basicnode] (v3) at (1,-1) {};
	\node [basicnode] (v4) at (-1.5,-2) {};
	\node [basicnode] (v5) at (-0.5,-2) {};
	\node [basicnode] (v6) at (0.5,-2) {};
	\node [basicnode] (v7) at (1.5,-2) {};

	\draw [blackedge] (v1) edge (v2);
	\draw [blackedge] (v1) edge (v3);
	\draw [blackedge] (v2) edge (v4);
	\draw [blackedge] (v2) edge (v5);
	\draw [blackedge] (v3) edge (v6);
	\draw [blackedge] (v3) edge (v7);
\end{scope}

\begin{scope}[shift={(9,-0.75)}]
 	\node [basicnode] (v1) at (0,0) {};
	\node [basicnode] (v2) at (-1,-1) {};
	\node [basicnode] (v3) at (1,-1) {};
	\node [basicnode] (v4) at (-1.5,-2) {};
	\node [basicnode] (v5) at (-0.5,-2) {};
	\node [basicnode] (v6) at (0.5,-2) {};
	\node [basicnode] (v7) at (1.5,-2) {};

	\draw [blackedge] (v2) edge (v5);
	\draw [blackedge] (v3) edge (v6);
	\draw [blackedge] (v4) edge (v5);
	\draw [blackedge] (v5) edge (v6);
	\draw [blackedge] (v6) edge (v7);
\end{scope}

\begin{scope}[shift={(14,-0.75)}]
 	\node [basicnode] (v1) at (0,0) {};
	\node [basicnode] (v2) at (-1,-1) {};
	\node [basicnode] (v3) at (1,-1) {};
	\node [basicnode] (v4) at (-1.5,-2) {};
	\node [basicnode] (v5) at (-0.5,-2) {};
	\node [basicnode] (v6) at (0.5,-2) {};
	\node [basicnode] (v7) at (1.5,-2) {};

	\draw [blackedge] (v1) edge (v2);
	\draw [blackedge] (v1) edge (v3);
	\draw [blackedge] (v2) edge (v4);
	\draw [blackedge] (v3) edge (v7);
	\draw [blackedge] (v4) edge (v5);
	\draw [blackedge] (v5) edge (v6);
	\draw [blackedge] (v6) edge (v7);
\end{scope}

\begin{scope}[shift={(0,-4)}]
	\node [blanknode] (h1) at (0, 0) {};
	\node [blanknode] (h2) at (18, 0) {};
	\draw [blackedge] (h1) edge (h2);
	
	\node [blanknode] (h3) at (9, 0.18) {};
	\node [blanknode] (h4) at (9, -7.93) {};
	\draw [blackedge] (h3) edge (h4);
	
	\node [blanknode] (h5) at (3.82, -7.75) {};
	\node [blanknode] (h6) at (9.18, -7.75) {};
	\draw [blackedge] (h5) edge (h6);
	
	\node [blanknode] (h5) at (4, -7.57) {};
	\node [blanknode] (h6) at (4, -11.43) {};
	\draw [blackedge] (h5) edge (h6);
	
	\node [blanknode] (h7) at (0, -11.25) {};
	\node [blanknode] (h8) at (18, -11.25) {};
	\draw [blackedge] (h7) edge (h8);
	
\end{scope}


\begin{scope}[shift={(2.5,-5.5)}]
 	\node [basicnode] (v1) at (0,0) {};
	\node [basicnode] (v2) at (-1,-1) {};
	\node [basicnode] (v3) at (1,-1) {};
	\node [basicnode] (v4) at (-1.5,-2) {};
	\node [basicnode] (v5) at (-0.5,-2) {};
	\node [basicnode] (v6) at (0.5,-2) {};
	\node [basicnode] (v7) at (1.5,-2) {};

	\draw [blackedge] (v1) edge (v2);
	\draw [blackedge] (v1) edge (v3);
	\draw [blackedge] (v2) edge (v4);
	\draw [rededge] (v2) edge (v5);
	\draw [rededge] (v3) edge (v6);
	\draw [blackedge] (v3) edge (v7);
	\draw [blueedge] (v4) edge (v5);
	\draw [blueedge] (v6) edge (v7);
	\draw [blueedge] (v5) edge (v6);
\end{scope}

\begin{scope}[shift={(6.5,-5.5)}]
 	\node [basicnode] (v1) at (0,0) {};
	\node [basicnode] (v2) at (-1,-1) {};
	\node [basicnode] (v3) at (1,-1) {};
	\node [basicnode] (v4) at (-1.5,-2) {};
	\node [basicnode] (v5) at (-0.5,-2) {};
	\node [basicnode] (v6) at (0.5,-2) {};
	\node [basicnode] (v7) at (1.5,-2) {};

	\draw [blackedge] (v1) edge (v2);
	\draw [blackedge] (v1) edge (v3);
	\draw [rededge] (v2) edge (v4);
	\draw [blackedge] (v2) edge (v5);
	\draw [blackedge] (v3) edge (v6);
	\draw [rededge] (v3) edge (v7);
	\draw [blueedge] (v4) edge (v5);
	\draw [blueedge] (v6) edge (v7);
	\draw [blueedge] (v4) edge[bend right=30] (v7);
\end{scope}

\begin{scope}[shift={(2.5,-9)}]
 	\node [basicnode] (v1) at (0,0) {};
	\node [basicnode] (v2) at (-1,-1) {};
	\node [basicnode] (v3) at (1,-1) {};
	\node [basicnode] (v4) at (-1.5,-2) {};
	\node [basicnode] (v5) at (-0.5,-2) {};
	\node [basicnode] (v6) at (0.5,-2) {};
	\node [basicnode] (v7) at (1.5,-2) {};

	\draw [blackedge] (v1) edge (v2);
	\draw [blackedge] (v1) edge (v3);
	\draw [rededge] (v2) edge (v4);
	\draw [blackedge] (v2) edge (v5);
	\draw [rededge] (v3) edge (v6);
	\draw [blackedge] (v3) edge (v7);
	\draw [blueedge] (v4) edge (v5);
	\draw [blueedge] (v6) edge (v7);
	\draw [blueedge] (v4) edge[bend right=35] (v6);
\end{scope}

\begin{scope}[shift={(6.5,-9)}]
 	\node [basicnode] (v1) at (0,0) {};
	\node [basicnode] (v2) at (-1,-1) {};
	\node [basicnode] (v3) at (1,-1) {};
	\node [basicnode] (v4) at (-1.5,-2) {};
	\node [basicnode] (v5) at (-0.5,-2) {};
	\node [basicnode] (v6) at (0.5,-2) {};
	\node [basicnode] (v7) at (1.5,-2) {};

	\draw [blackedge] (v1) edge (v2);
	\draw [blackedge] (v1) edge (v3);
	\draw [blackedge] (v2) edge (v4);
	\draw [rededge] (v2) edge (v5);
	\draw [blackedge] (v3) edge (v6);
	\draw [rededge] (v3) edge (v7);
	\draw [blueedge] (v4) edge (v5);
	\draw [blueedge] (v6) edge (v7);
	\draw [blueedge] (v5) edge[bend right=35] (v7);
\end{scope}

\begin{scope}[shift={(11.5,-5.5)}]
 	\node [basicnode] (v1) at (0,0) {};
	\node [basicnode] (v2) at (-1,-1) {};
	\node [basicnode] (v3) at (1,-1) {};
	\node [basicnode] (v4) at (-1.5,-2) {};
	\node [basicnode] (v5) at (-0.5,-2) {};
	\node [basicnode] (v6) at (0.5,-2) {};
	\node [basicnode] (v7) at (1.5,-2) {};

	\draw [blackedge] (v1) edge (v2);
	\draw [blackedge] (v1) edge (v3);
	\draw [rededge] (v2) edge (v4);
	\draw [blackedge] (v3) edge (v7);
	\draw [rededge] (v4) edge (v5);
	\draw [rededge] (v5) edge (v6);
	\draw [blackedge] (v6) edge (v7);
	
	\draw [blueedge] (v4) edge[bend right=20] (v1);
	\draw [blueedge] (v5) edge[bend left=35] (v7);
\end{scope}

\begin{scope}[shift={(15.5,-5.5)}]
 	\node [basicnode] (v1) at (0,0) {};
	\node [basicnode] (v2) at (-1,-1) {};
	\node [basicnode] (v3) at (1,-1) {};
	\node [basicnode] (v4) at (-1.5,-2) {};
	\node [basicnode] (v5) at (-0.5,-2) {};
	\node [basicnode] (v6) at (0.5,-2) {};
	\node [basicnode] (v7) at (1.5,-2) {};

	\draw [blackedge] (v1) edge (v2);
	\draw [blackedge] (v1) edge (v3);
	\draw [blackedge] (v2) edge (v4);
	\draw [rededge] (v3) edge (v7);
	\draw [blackedge] (v4) edge (v5);
	\draw [rededge] (v5) edge (v6);
	\draw [rededge] (v6) edge (v7);
	
	\draw [blueedge] (v6) edge[bend right=35] (v4);
	\draw [blueedge] (v7) edge[bend left=20] (v1);
\end{scope}

\begin{scope}[shift={(11.5,-9)}]
 	\node [basicnode] (v1) at (0,0) {};
	\node [basicnode] (v2) at (-1,-1) {};
	\node [basicnode] (v3) at (1,-1) {};
	\node [basicnode] (v4) at (-1.5,-2) {};
	\node [basicnode] (v5) at (-0.5,-2) {};
	\node [basicnode] (v6) at (0.5,-2) {};
	\node [basicnode] (v7) at (1.5,-2) {};

	\draw [rededge] (v1) edge (v2);
	\draw [blackedge] (v1) edge (v3);
	\draw [rededge] (v2) edge (v4);
	\draw [blackedge] (v3) edge (v7);
	\draw [rededge] (v4) edge (v5);
	\draw [blackedge] (v5) edge (v6);
	\draw [blackedge] (v6) edge (v7);
	
	\draw [blueedge] (v2) edge (v3);
	\draw [blueedge] (v4) edge[bend left=35] (v6);
\end{scope}

\begin{scope}[shift={(15.5,-9)}]
 	\node [basicnode] (v1) at (0,0) {};
	\node [basicnode] (v2) at (-1,-1) {};
	\node [basicnode] (v3) at (1,-1) {};
	\node [basicnode] (v4) at (-1.5,-2) {};
	\node [basicnode] (v5) at (-0.5,-2) {};
	\node [basicnode] (v6) at (0.5,-2) {};
	\node [basicnode] (v7) at (1.5,-2) {};

	\draw [blackedge] (v1) edge (v2);
	\draw [rededge] (v1) edge (v3);
	\draw [blackedge] (v2) edge (v4);
	\draw [rededge] (v3) edge (v7);
	\draw [blackedge] (v4) edge (v5);
	\draw [blackedge] (v5) edge (v6);
	\draw [rededge] (v6) edge (v7);
	
	\draw [blueedge] (v7) edge[bend right=35] (v5);
	\draw [blueedge] (v3) edge (v2);
\end{scope}

\begin{scope}[shift={(7.5,-12.5)}]
 	\node [basicnode] (v1) at (0,0) {};
	\node [basicnode] (v2) at (-1,-1) {};
	\node [basicnode] (v3) at (1,-1) {};
	\node [basicnode] (v4) at (-1.5,-2) {};
	\node [basicnode] (v5) at (-0.5,-2) {};
	\node [basicnode] (v6) at (0.5,-2) {};
	\node [basicnode] (v7) at (1.5,-2) {};

	\draw [blackedge] (v1) edge (v2);
	\draw [blackedge] (v1) edge (v3);
	\draw [blackedge] (v2) edge (v4);
	\draw [blackedge] (v3) edge (v7);
	\draw [rededge] (v4) edge (v5);
	\draw [rededge] (v5) edge (v6);
	\draw [rededge] (v6) edge (v7);
	
	\draw [blueedge] (v2) edge (v5);
	\draw [blueedge] (v3) edge (v6);
\end{scope}

\begin{scope}[shift={(11.5,-12.5)}]
 	\node [basicnode] (v1) at (0,0) {};
	\node [basicnode] (v2) at (-1,-1) {};
	\node [basicnode] (v3) at (1,-1) {};
	\node [basicnode] (v4) at (-1.5,-2) {};
	\node [basicnode] (v5) at (-0.5,-2) {};
	\node [basicnode] (v6) at (0.5,-2) {};
	\node [basicnode] (v7) at (1.5,-2) {};

	\draw [rededge] (v1) edge (v2);
	\draw [rededge] (v1) edge (v3);
	\draw [rededge] (v2) edge (v4);
	\draw [blackedge] (v3) edge (v7);
	\draw [blackedge] (v4) edge (v5);
	\draw [blackedge] (v5) edge (v6);
	\draw [blackedge] (v6) edge (v7);
	
	\draw [blueedge] (v1) edge[bend right=20] (v7);
	\draw [blueedge] (v2) edge (v5);
\end{scope}

\begin{scope}[shift={(15.5,-12.5)}]
 	\node [basicnode] (v1) at (0,0) {};
	\node [basicnode] (v2) at (-1,-1) {};
	\node [basicnode] (v3) at (1,-1) {};
	\node [basicnode] (v4) at (-1.5,-2) {};
	\node [basicnode] (v5) at (-0.5,-2) {};
	\node [basicnode] (v6) at (0.5,-2) {};
	\node [basicnode] (v7) at (1.5,-2) {};

	\draw [rededge] (v1) edge (v2);
	\draw [rededge] (v1) edge (v3);
	\draw [blackedge] (v2) edge (v4);
	\draw [rededge] (v3) edge (v7);
	\draw [blackedge] (v4) edge (v5);
	\draw [blackedge] (v5) edge (v6);
	\draw [blackedge] (v6) edge (v7);
	
	\draw [blueedge] (v3) edge (v6);
	\draw [blueedge] (v1) edge[bend left=20] (v4);
\end{scope}

\end{tikzpicture}

%% file: figures/edge_set_orbit_example.tex
\begin{tikzpicture}[scale=0.75]

\node [textnode] at (-4.5,0) {Example Graph $G$};
\node [textnode] at (4.5,0) {Automorphism Orbit of};
\node [textnode] at (4.5,-0.75) {Edge Set $\{(1, 5), (5, 6)\}$ in $G$};

\node [textnode] at (0,-6) {AO$_G\left( \{(1, 5), (5, 6)\} \right) = \{\textcolor{blue}{\textbf{\{(1, 5), (5, 6)\}}},\ \  \textcolor{green!60!black}{\textbf{\{(6, 7), (7, 3)\}}},$};
\node [textnode] at (2.85, -6.75) {$\textcolor{red!80!brown}{\textbf{\{(2, 5), (5, 6)\}}},\ \  \textcolor{gray!80!black}{\textbf{\{(6, 7), (7, 4)\}}}\}$};

\node [textnode] at (0,-7.75) {Stab$_G\left( \{(1, 5), (5, 6)\} \right) = \{(1 \rightarrow 1, 2 \rightarrow 2, 3 \rightarrow 3, 4 \rightarrow 4, 5 \rightarrow 5, 6 \rightarrow 6, 7 \rightarrow 7),$};
\node [textnode] at (2.95, -8.5) {$(1 \rightarrow 1, 2 \rightarrow 2, 3 \rightarrow 4, 4 \rightarrow 3, 5 \rightarrow 5, 6 \rightarrow 6, 7 \rightarrow 7)\}$};

\begin{scope}[shift={(-4.5,-4)}]
\node [basicnode, label=6] (v6) at (0,1) {};
\node [basicnode, label=7] (v7) at (1.5,1) {};
\node [basicnode, label=5] (v5) at (-1.5,1) {};
\node [basicnode, label=1] (v1) at (-3,2) {};
\node [basicnode, label=2] (v2) at (-3,0) {};
\node [basicnode, label=3] (v3) at (3,2) {};
\node [basicnode, label=4] (v4) at (3,0) {};

\draw [basicedge] (v2) edge (v5);
\draw [basicedge] (v1) edge (v5);
\draw [basicedge] (v5) edge (v6);
\draw [basicedge] (v6) edge (v7);
\draw [basicedge] (v7) edge (v3);
\draw [basicedge] (v7) edge (v4);
\end{scope}

\begin{scope}[shift={(4.5,-4)}]
\node [basicnode, label=\textcolor{lightgray}{6}] (v6) at (0,1) {};
\node [basicnode, label=\textcolor{lightgray}{7}] (v7) at (1.5,1) {};
\node [basicnode, label=\textcolor{lightgray}{5}] (v5) at (-1.5,1) {};
\node [basicnode, label=\textcolor{lightgray}{1}] (v1) at (-3,2) {};
\node [basicnode, label=\textcolor{lightgray}{2}] (v2) at (-3,0) {};
\node [basicnode, label=\textcolor{lightgray}{3}] (v3) at (3,2) {};
\node [basicnode, label=\textcolor{lightgray}{4}] (v4) at (3,0) {};

\draw [basicedge, draw=blue] (v1) edge["a", bend left=15] (v5); 
\draw [basicedge, draw=blue] (v5) edge["a", bend left=20] (v6); 

\draw [basicedge, draw=red!80!brown] (v5) edge["b", bend left=15] (v2); 
\draw [basicedge, draw=red!80!brown] (v6) edge["b", bend left=20] (v5); 

\draw [basicedge, draw=green!60!black] (v6) edge["c", bend left=20] (v7); 
\draw [basicedge, draw=green!60!black] (v7) edge["c", bend left=15] (v3); 

\draw [basicedge, draw=gray!80!black] (v7) edge["d", bend left=20] (v6); 
\draw [basicedge, draw=gray!80!black] (v4) edge["d", bend left=15] (v7); 
\end{scope}


\begin{scope}[shift={(9,0)}]
\node [blanknode] (v9) at (-9,-1) {};
\node [blanknode] (v8) at (-9,-5) {};
\draw [basicedge] (v8) edge (v9);
\end{scope}

\end{tikzpicture}

%% file: appendix.tex
\section*{Appendix -- Directed Example}
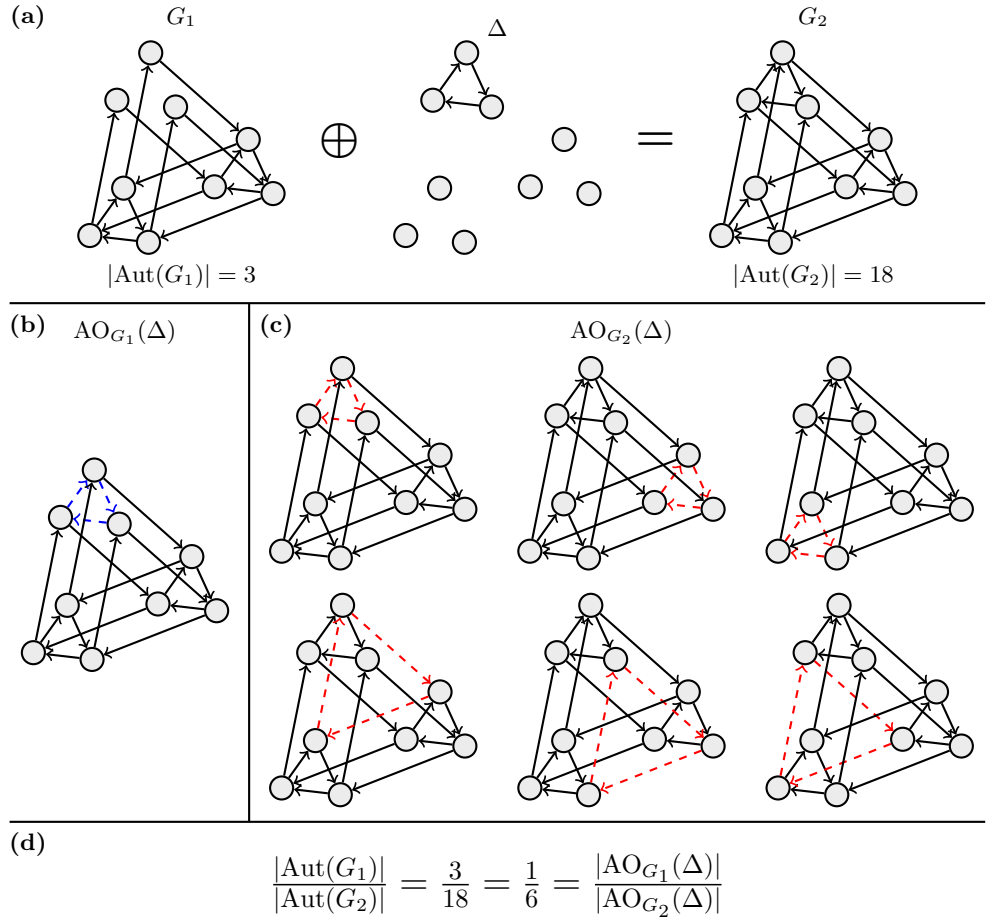
\begin{figure}[h]
\centering
\input{figures/main_result_03}
\caption{\textbf{Directed Example of Main Result} -- \textbf{(a)} shows two graphs $G_1$ and $G_2$, along with the difference between them $\Delta$. \textbf{(b) and (c)} show the automorphism orbits of $\Delta$ in $G_1$ and $G_2$ respectively. Red lines indicate edges (to be deleted). Blue lines indicate non-edges (to be added).
Finally, \textbf{(d)} shows our formula applied to this example.\\
It might not be immediately obvious that $G_2$ has 18 automorphisms; the first 9 are easier to note. However, as the depiction of $\text{AO}_{G_2}(\Delta)$ shows, each little triangle in $G_2$ is structurally equivalent to each big triangle. This graph is similar to the Peterson graph in that the equivalency of some cycles is not immediately apparent upon visual inspection.}\label{fig:main_example_03}
\end{figure}

%% file: figures/main_result_03.tex
\begin{tikzpicture}[scale=0.6]


\node [textnode] at (-1.4, 0.25) {\textbf{(a)}};
\node [textnode] at (-1.4, -6.6) {\textbf{(b)}};
\node [textnode] at (4.1, -6.6) {\textbf{(c)}};
\node [textnode] at (-1.4, -18.1) {\textbf{(d)}};

\node [textnode] at (2,0.25) {$G_1$};
\node [textnode] at (5.5,-2.5) {\LARGE $\bm{\oplus}$};
\node [textnode] at (9,0) {$\Delta$};
\node [textnode] at (12.5,-2.5) {\LARGE $\bm{=}$};
\node [textnode] at (16,0.25) {$G_2$};

\begin{scope}[shift={(0,-5.5)}]
	\node [textnode] at (2, 0) {$|\text{Aut}(G_1)| = 3$};
	\node [textnode] at (16, 0) {$|\text{Aut}(G_2)| = 18$};
\end{scope}

\begin{scope}[shift={(0,-6.75)}]
	\node [textnode] at (0.75, 0) {$\text{AO}_{G_1}(\Delta)$};
	\node [textnode] at (11.75, 0) {$\text{AO}_{G_2}(\Delta)$};
\end{scope}

\node [textnode] at (9,-19) {\Large $\frac{|\text{Aut}(G_1)|}{|\text{Aut}(G_2)|} = \frac{3}{18} = \frac{1}{6} = \frac{|\text{AO}_{G_1}(\Delta)|}{|\text{AO}_{G_2}(\Delta)|}$};

\begin{scope}[shift={(-0.03,-4.75)}]

\begin{scope}[shift={(0, 0)}]
 	\node [basicnode] (v1) at (0, 0.15) {};
	\node [basicnode] (v2) at (1.3, 0) {};
	\node [basicnode] (v3) at (0.75, 1.2) {};
\end{scope}

\begin{scope}[shift={(0.6, 3)}]
 	\node [basicnode] (v4) at (0, 0.15) {};
	\node [basicnode] (v5) at (1.3, 0) {};
	\node [basicnode] (v6) at (0.75, 1.2) {};
\end{scope}

\begin{scope}[shift={(2.76, 1.08)}]
 	\node [basicnode] (v7) at (0, 0.15) {};
	\node [basicnode] (v8) at (1.3, 0) {};
	\node [basicnode] (v9) at (0.75, 1.2) {};
\end{scope}

	\draw [blackdiredge] (v1) edge (v3);
	\draw [blackdiredge] (v3) edge (v2);
	\draw [blackdiredge] (v2) edge (v1);
	\draw [blackdiredge] (v7) edge (v9);
	\draw [blackdiredge] (v9) edge (v8);
	\draw [blackdiredge] (v8) edge (v7);
	
	\draw [blackdiredge] (v1) edge (v4);
	\draw [blackdiredge] (v4) edge (v7);
	\draw [blackdiredge] (v7) edge (v1);
	\draw [blackdiredge] (v2) edge (v5);
	\draw [blackdiredge] (v5) edge (v8);
	\draw [blackdiredge] (v8) edge (v2);
	\draw [blackdiredge] (v3) edge (v6);
	\draw [blackdiredge] (v6) edge (v9);
	\draw [blackdiredge] (v9) edge (v3);
\end{scope}

\begin{scope}[shift={(6.97,-4.75)}]

\begin{scope}[shift={(0, 0)}]
 	\node [basicnode] (v1) at (0, 0.15) {};
	\node [basicnode] (v2) at (1.3, 0) {};
	\node [basicnode] (v3) at (0.75, 1.2) {};
\end{scope}

\begin{scope}[shift={(0.6, 3)}]
 	\node [basicnode] (v4) at (0, 0.15) {};
	\node [basicnode] (v5) at (1.3, 0) {};
	\node [basicnode] (v6) at (0.75, 1.2) {};
\end{scope}

\begin{scope}[shift={(2.76, 1.08)}]
 	\node [basicnode] (v7) at (0, 0.15) {};
	\node [basicnode] (v8) at (1.3, 0) {};
	\node [basicnode] (v9) at (0.75, 1.2) {};
\end{scope}

	\draw [blackdiredge] (v4) edge (v6);
	\draw [blackdiredge] (v6) edge (v5);
	\draw [blackdiredge] (v5) edge (v4);
\end{scope}

\begin{scope}[shift={(13.97,-4.75)}]

\begin{scope}[shift={(0, 0)}]
 	\node [basicnode] (v1) at (0, 0.15) {};
	\node [basicnode] (v2) at (1.3, 0) {};
	\node [basicnode] (v3) at (0.75, 1.2) {};
\end{scope}

\begin{scope}[shift={(0.6, 3)}]
 	\node [basicnode] (v4) at (0, 0.15) {};
	\node [basicnode] (v5) at (1.3, 0) {};
	\node [basicnode] (v6) at (0.75, 1.2) {};
\end{scope}

\begin{scope}[shift={(2.76, 1.08)}]
 	\node [basicnode] (v7) at (0, 0.15) {};
	\node [basicnode] (v8) at (1.3, 0) {};
	\node [basicnode] (v9) at (0.75, 1.2) {};
\end{scope}

	\draw [blackdiredge] (v1) edge (v3);
	\draw [blackdiredge] (v3) edge (v2);
	\draw [blackdiredge] (v2) edge (v1);
	\draw [blackdiredge] (v4) edge (v6);
	\draw [blackdiredge] (v6) edge (v5);
	\draw [blackdiredge] (v5) edge (v4);
	\draw [blackdiredge] (v7) edge (v9);
	\draw [blackdiredge] (v9) edge (v8);
	\draw [blackdiredge] (v8) edge (v7);
	
	\draw [blackdiredge] (v1) edge (v4);
	\draw [blackdiredge] (v4) edge (v7);
	\draw [blackdiredge] (v7) edge (v1);
	\draw [blackdiredge] (v2) edge (v5);
	\draw [blackdiredge] (v5) edge (v8);
	\draw [blackdiredge] (v8) edge (v2);
	\draw [blackdiredge] (v3) edge (v6);
	\draw [blackdiredge] (v6) edge (v9);
	\draw [blackdiredge] (v9) edge (v3);
\end{scope}

\begin{scope}[shift={(0,-6.1)}]
	\node [blanknode] (h1) at (-2, 0) {};
	\node [blanknode] (h2) at (20, 0) {};
	\draw [blackedge] (h1) edge (h2);
	
	\node [blanknode] (h5) at (3.5, 0.18) {};
	\node [blanknode] (h6) at (3.5, -11.68) {};
	\draw [blackedge] (h5) edge (h6);
	
	\node [blanknode] (h7) at (-2, -11.5) {};
	\node [blanknode] (h8) at (20, -11.5) {};
	\draw [blackedge] (h7) edge (h8);
	
\end{scope}


\begin{scope}[shift={(-1.28,-14)}]

\begin{scope}[shift={(0, 0)}]
 	\node [basicnode] (v1) at (0, 0.15) {};
	\node [basicnode] (v2) at (1.3, 0) {};
	\node [basicnode] (v3) at (0.75, 1.2) {};
\end{scope}

\begin{scope}[shift={(0.6, 3)}]
 	\node [basicnode] (v4) at (0, 0.15) {};
	\node [basicnode] (v5) at (1.3, 0) {};
	\node [basicnode] (v6) at (0.75, 1.2) {};
\end{scope}

\begin{scope}[shift={(2.76, 1.08)}]
 	\node [basicnode] (v7) at (0, 0.15) {};
	\node [basicnode] (v8) at (1.3, 0) {};
	\node [basicnode] (v9) at (0.75, 1.2) {};
\end{scope}

	\draw [blackdiredge] (v1) edge (v3);
	\draw [blackdiredge] (v3) edge (v2);
	\draw [blackdiredge] (v2) edge (v1);
	\draw [bluediredge] (v4) edge (v6);
	\draw [bluediredge] (v6) edge (v5);
	\draw [bluediredge] (v5) edge (v4);
	\draw [blackdiredge] (v7) edge (v9);
	\draw [blackdiredge] (v9) edge (v8);
	\draw [blackdiredge] (v8) edge (v7);
	
	\draw [blackdiredge] (v1) edge (v4);
	\draw [blackdiredge] (v4) edge (v7);
	\draw [blackdiredge] (v7) edge (v1);
	\draw [blackdiredge] (v2) edge (v5);
	\draw [blackdiredge] (v5) edge (v8);
	\draw [blackdiredge] (v8) edge (v2);
	\draw [blackdiredge] (v3) edge (v6);
	\draw [blackdiredge] (v6) edge (v9);
	\draw [blackdiredge] (v9) edge (v3);
\end{scope}

\begin{scope}[shift={(4.22,-11.75)}]

\begin{scope}[shift={(0, 0)}]
 	\node [basicnode] (v1) at (0, 0.15) {};
	\node [basicnode] (v2) at (1.3, 0) {};
	\node [basicnode] (v3) at (0.75, 1.2) {};
\end{scope}

\begin{scope}[shift={(0.6, 3)}]
 	\node [basicnode] (v4) at (0, 0.15) {};
	\node [basicnode] (v5) at (1.3, 0) {};
	\node [basicnode] (v6) at (0.75, 1.2) {};
\end{scope}

\begin{scope}[shift={(2.76, 1.08)}]
 	\node [basicnode] (v7) at (0, 0.15) {};
	\node [basicnode] (v8) at (1.3, 0) {};
	\node [basicnode] (v9) at (0.75, 1.2) {};
\end{scope}

	\draw [blackdiredge] (v1) edge (v3);
	\draw [blackdiredge] (v3) edge (v2);
	\draw [blackdiredge] (v2) edge (v1);
	\draw [reddiredge] (v4) edge (v6);
	\draw [reddiredge] (v6) edge (v5);
	\draw [reddiredge] (v5) edge (v4);
	\draw [blackdiredge] (v7) edge (v9);
	\draw [blackdiredge] (v9) edge (v8);
	\draw [blackdiredge] (v8) edge (v7);
	
	\draw [blackdiredge] (v1) edge (v4);
	\draw [blackdiredge] (v4) edge (v7);
	\draw [blackdiredge] (v7) edge (v1);
	\draw [blackdiredge] (v2) edge (v5);
	\draw [blackdiredge] (v5) edge (v8);
	\draw [blackdiredge] (v8) edge (v2);
	\draw [blackdiredge] (v3) edge (v6);
	\draw [blackdiredge] (v6) edge (v9);
	\draw [blackdiredge] (v9) edge (v3);
\end{scope}

\begin{scope}[shift={(9.72,-11.75)}]

\begin{scope}[shift={(0, 0)}]
 	\node [basicnode] (v1) at (0, 0.15) {};
	\node [basicnode] (v2) at (1.3, 0) {};
	\node [basicnode] (v3) at (0.75, 1.2) {};
\end{scope}

\begin{scope}[shift={(0.6, 3)}]
 	\node [basicnode] (v4) at (0, 0.15) {};
	\node [basicnode] (v5) at (1.3, 0) {};
	\node [basicnode] (v6) at (0.75, 1.2) {};
\end{scope}

\begin{scope}[shift={(2.76, 1.08)}]
 	\node [basicnode] (v7) at (0, 0.15) {};
	\node [basicnode] (v8) at (1.3, 0) {};
	\node [basicnode] (v9) at (0.75, 1.2) {};
\end{scope}

	\draw [blackdiredge] (v1) edge (v3);
	\draw [blackdiredge] (v3) edge (v2);
	\draw [blackdiredge] (v2) edge (v1);
	\draw [blackdiredge] (v4) edge (v6);
	\draw [blackdiredge] (v6) edge (v5);
	\draw [blackdiredge] (v5) edge (v4);
	\draw [reddiredge] (v7) edge (v9);
	\draw [reddiredge] (v9) edge (v8);
	\draw [reddiredge] (v8) edge (v7);
	
	\draw [blackdiredge] (v1) edge (v4);
	\draw [blackdiredge] (v4) edge (v7);
	\draw [blackdiredge] (v7) edge (v1);
	\draw [blackdiredge] (v2) edge (v5);
	\draw [blackdiredge] (v5) edge (v8);
	\draw [blackdiredge] (v8) edge (v2);
	\draw [blackdiredge] (v3) edge (v6);
	\draw [blackdiredge] (v6) edge (v9);
	\draw [blackdiredge] (v9) edge (v3);
\end{scope}

\begin{scope}[shift={(15.22,-11.75)}]

\begin{scope}[shift={(0, 0)}]
 	\node [basicnode] (v1) at (0, 0.15) {};
	\node [basicnode] (v2) at (1.3, 0) {};
	\node [basicnode] (v3) at (0.75, 1.2) {};
\end{scope}

\begin{scope}[shift={(0.6, 3)}]
 	\node [basicnode] (v4) at (0, 0.15) {};
	\node [basicnode] (v5) at (1.3, 0) {};
	\node [basicnode] (v6) at (0.75, 1.2) {};
\end{scope}

\begin{scope}[shift={(2.76, 1.08)}]
 	\node [basicnode] (v7) at (0, 0.15) {};
	\node [basicnode] (v8) at (1.3, 0) {};
	\node [basicnode] (v9) at (0.75, 1.2) {};
\end{scope}

	\draw [reddiredge] (v1) edge (v3);
	\draw [reddiredge] (v3) edge (v2);
	\draw [reddiredge] (v2) edge (v1);
	\draw [blackdiredge] (v4) edge (v6);
	\draw [blackdiredge] (v6) edge (v5);
	\draw [blackdiredge] (v5) edge (v4);
	\draw [blackdiredge] (v7) edge (v9);
	\draw [blackdiredge] (v9) edge (v8);
	\draw [blackdiredge] (v8) edge (v7);
	
	\draw [blackdiredge] (v1) edge (v4);
	\draw [blackdiredge] (v4) edge (v7);
	\draw [blackdiredge] (v7) edge (v1);
	\draw [blackdiredge] (v2) edge (v5);
	\draw [blackdiredge] (v5) edge (v8);
	\draw [blackdiredge] (v8) edge (v2);
	\draw [blackdiredge] (v3) edge (v6);
	\draw [blackdiredge] (v6) edge (v9);
	\draw [blackdiredge] (v9) edge (v3);
\end{scope}

\begin{scope}[shift={(4.22,-17)}]

\begin{scope}[shift={(0, 0)}]
 	\node [basicnode] (v1) at (0, 0.15) {};
	\node [basicnode] (v2) at (1.3, 0) {};
	\node [basicnode] (v3) at (0.75, 1.2) {};
\end{scope}

\begin{scope}[shift={(0.6, 3)}]
 	\node [basicnode] (v4) at (0, 0.15) {};
	\node [basicnode] (v5) at (1.3, 0) {};
	\node [basicnode] (v6) at (0.75, 1.2) {};
\end{scope}

\begin{scope}[shift={(2.76, 1.08)}]
 	\node [basicnode] (v7) at (0, 0.15) {};
	\node [basicnode] (v8) at (1.3, 0) {};
	\node [basicnode] (v9) at (0.75, 1.2) {};
\end{scope}

	\draw [blackdiredge] (v1) edge (v3);
	\draw [blackdiredge] (v3) edge (v2);
	\draw [blackdiredge] (v2) edge (v1);
	\draw [blackdiredge] (v4) edge (v6);
	\draw [blackdiredge] (v6) edge (v5);
	\draw [blackdiredge] (v5) edge (v4);
	\draw [blackdiredge] (v7) edge (v9);
	\draw [blackdiredge] (v9) edge (v8);
	\draw [blackdiredge] (v8) edge (v7);
	
	\draw [blackdiredge] (v1) edge (v4);
	\draw [blackdiredge] (v4) edge (v7);
	\draw [blackdiredge] (v7) edge (v1);
	\draw [blackdiredge] (v2) edge (v5);
	\draw [blackdiredge] (v5) edge (v8);
	\draw [blackdiredge] (v8) edge (v2);
	\draw [reddiredge] (v3) edge (v6);
	\draw [reddiredge] (v6) edge (v9);
	\draw [reddiredge] (v9) edge (v3);
\end{scope}

\begin{scope}[shift={(9.72,-17)}]

\begin{scope}[shift={(0, 0)}]
 	\node [basicnode] (v1) at (0, 0.15) {};
	\node [basicnode] (v2) at (1.3, 0) {};
	\node [basicnode] (v3) at (0.75, 1.2) {};
\end{scope}

\begin{scope}[shift={(0.6, 3)}]
 	\node [basicnode] (v4) at (0, 0.15) {};
	\node [basicnode] (v5) at (1.3, 0) {};
	\node [basicnode] (v6) at (0.75, 1.2) {};
\end{scope}

\begin{scope}[shift={(2.76, 1.08)}]
 	\node [basicnode] (v7) at (0, 0.15) {};
	\node [basicnode] (v8) at (1.3, 0) {};
	\node [basicnode] (v9) at (0.75, 1.2) {};
\end{scope}

	\draw [blackdiredge] (v1) edge (v3);
	\draw [blackdiredge] (v3) edge (v2);
	\draw [blackdiredge] (v2) edge (v1);
	\draw [blackdiredge] (v4) edge (v6);
	\draw [blackdiredge] (v6) edge (v5);
	\draw [blackdiredge] (v5) edge (v4);
	\draw [blackdiredge] (v7) edge (v9);
	\draw [blackdiredge] (v9) edge (v8);
	\draw [blackdiredge] (v8) edge (v7);
	
	\draw [blackdiredge] (v1) edge (v4);
	\draw [blackdiredge] (v4) edge (v7);
	\draw [blackdiredge] (v7) edge (v1);
	\draw [reddiredge] (v2) edge (v5);
	\draw [reddiredge] (v5) edge (v8);
	\draw [reddiredge] (v8) edge (v2);
	\draw [blackdiredge] (v3) edge (v6);
	\draw [blackdiredge] (v6) edge (v9);
	\draw [blackdiredge] (v9) edge (v3);
\end{scope}

\begin{scope}[shift={(15.22,-17)}]

\begin{scope}[shift={(0, 0)}]
 	\node [basicnode] (v1) at (0, 0.15) {};
	\node [basicnode] (v2) at (1.3, 0) {};
	\node [basicnode] (v3) at (0.75, 1.2) {};
\end{scope}

\begin{scope}[shift={(0.6, 3)}]
 	\node [basicnode] (v4) at (0, 0.15) {};
	\node [basicnode] (v5) at (1.3, 0) {};
	\node [basicnode] (v6) at (0.75, 1.2) {};
\end{scope}

\begin{scope}[shift={(2.76, 1.08)}]
 	\node [basicnode] (v7) at (0, 0.15) {};
	\node [basicnode] (v8) at (1.3, 0) {};
	\node [basicnode] (v9) at (0.75, 1.2) {};
\end{scope}

	\draw [blackdiredge] (v1) edge (v3);
	\draw [blackdiredge] (v3) edge (v2);
	\draw [blackdiredge] (v2) edge (v1);
	\draw [blackdiredge] (v4) edge (v6);
	\draw [blackdiredge] (v6) edge (v5);
	\draw [blackdiredge] (v5) edge (v4);
	\draw [blackdiredge] (v7) edge (v9);
	\draw [blackdiredge] (v9) edge (v8);
	\draw [blackdiredge] (v8) edge (v7);
	
	\draw [reddiredge] (v1) edge (v4);
	\draw [reddiredge] (v4) edge (v7);
	\draw [reddiredge] (v7) edge (v1);
	\draw [blackdiredge] (v2) edge (v5);
	\draw [blackdiredge] (v5) edge (v8);
	\draw [blackdiredge] (v8) edge (v2);
	\draw [blackdiredge] (v3) edge (v6);
	\draw [blackdiredge] (v6) edge (v9);
	\draw [blackdiredge] (v9) edge (v3);
\end{scope}

\end{tikzpicture}

%% file: arXiv_main.bbl
\begin{thebibliography}{10}

\bibitem{andersen1996set}
Lars~D{\o}vling Andersen, Songkang Ding, and Preben~Dahl Vestergaard.
\newblock On the set edge-reconstruction conjecture.
\newblock {\em JCMCC-Journal of Combinatorial Mathematics and Combinatorial
  Computing}, 20:3--9, 1996.

\bibitem{bollobas1990almost}
B{\'e}la Bollob{\'a}s.
\newblock Almost every graph has reconstruction number three.
\newblock {\em Journal of Graph Theory}, 14(1):1--4, 1990.

\bibitem{bondy1977graph}
J~Adrian Bondy and Robert~L Hemminger.
\newblock Graph reconstruction—a survey.
\newblock {\em Journal of Graph Theory}, 1(3):227--268, 1977.

\bibitem{greenwell1969reconstructing}
DL~Greenwell and RL~Hemminger.
\newblock Reconstructing graphs.
\newblock In {\em The many facets of graph theory}, pages 91--114. Springer,
  1969.

\bibitem{harary2014graphical}
Frank Harary and Edgar~M Palmer.
\newblock {\em Graphical enumeration}.
\newblock Elsevier, 2014.

\bibitem{harary1967number}
Frank Harary, Edgar~M Palmer, and Ronald~C Read.
\newblock The number of ways to label a structure.
\newblock {\em Psychometrika}, 32(2):155--156, 1967.

\bibitem{hartke2009automorphism}
Stephen~G Hartke, Hannah Kolb, Jared Nishikawa, and Derrick Stolee.
\newblock Automorphism groups of a graph and a vertex-deleted subgraph.
\newblock {\em Electronic Journal of Combinatorics}, 2010.

\bibitem{kocay1987family}
William~L Kocay.
\newblock A family of nonreconstructible hypergraphs.
\newblock {\em Journal of Combinatorial Theory, Series B}, 42(1):46--63, 1987.

\bibitem{mcmullen2006graph}
Brian McMullen.
\newblock Graph reconstruction numbers.
\newblock 2006.

\bibitem{Note1}
Here $\oplus $ is analogous to the XOR function.

\bibitem{rose2009action}
HE~Rose.
\newblock Action and the orbit--stabiliser theorem.
\newblock {\em A Course on Finite Groups}, pages 91--111, 2009.

\end{thebibliography}
